\documentclass[english,12pt]{article}
\usepackage[T1]{fontenc}
\usepackage[latin1]{inputenc}
\usepackage{amsmath}
\usepackage{amsthm}
\usepackage{amssymb}
\usepackage{graphicx}

\makeatletter
 \theoremstyle{plain}
 \newtheorem{thm}{Theorem}

 \numberwithin{equation}{section} %% Comment out for sequentially-numbered

 \numberwithin{figure}{section} %% Comment out for sequentially-numbered

 \theoremstyle{plain}

 \theoremstyle{remark}

 \theoremstyle{plain}  
 \theoremstyle{plain}

 \theoremstyle{plain}

 \theoremstyle{plain}

 \theoremstyle{plain}
 \theoremstyle{plain}
 \newtheorem{lemma}{Lemma} %%Delete [thm] to re-start numbering

 \theoremstyle{definition}

 \theoremstyle{definition}

\usepackage{babel}
\makeatother

\newcommand{\tmop}[1]{\mathrm{#1}}
\renewcommand{\dfrac}[0]{\displaystyle\frac}

\begin{document}

\newcommand{\Ub}{\overline{U}}
\newcommand{\Pb}{\overline{P}}

\title{Nonexistence of Local Self-Similar Blow-up for the 3D Incompressible 
Navier-Stokes Equations}
\author{Thomas Y. Hou\thanks{Applied and Comput. Math, Caltech, Pasadena,
CA 91125. Email: hou@acm.caltech.edu.} 
and Ruo Li\thanks{Applied and Comput. Math., Caltech, Pasadena, CA 91125, and
LMAM\&School of Mathematical Sciences, Peking University, Beijing, China.
Email: rli@acm.caltech.edu.}} 

\maketitle

\abstract{We prove the nonexistence of local self-similar solutions of
the three dimensional incompressible Navier-Stokes equations.
The local self-similar solutions we consider here are different
from the global self-similar solutions. The self-similar scaling 
is only valid in an inner core region which shrinks to a point 
dynamically as the time, $t$, approaches the singularity time, $T$. 
The solution outside the inner core region is assumed to be regular. 
Under the assumption that the local self-similar velocity profile 
converges to a limiting profile as $t \rightarrow T$ in $L^p$ for some
$p \in (3,\infty )$, we prove that such local self-similar blow-up 
is not possible for any finite time.}

\section{Introduction.}

In this paper, we study the 3D incompressible Navier-Stokes equations 
with unit viscosity and zero external force:
\begin{equation}
  \left\{ \begin{array}{l}
    u_t + (u \cdot \nabla ) u = - \nabla p + \Delta u,\\
    \nabla \cdot u = 0,\\
    u|_{t = 0} = u_0 ( x ).
  \end{array} \label{nse} \right.
\end{equation}
We assume that the initial condition $u_0$ is divergence free
and $u_0 \in L^2(\mathbb{R}^3)\cap L^p(\mathbb{R}^3)$ for some 
$p\in (3,\infty)$.

Many physicists and mathematicians have made a great deal of effort 
in understanding the physical as well as the mathematical properties 
of the 3D incompressible Navier-Stokes equations. One of the long 
standing open questions is whether the solution of the 3D Navier-Stokes 
equations can develop a finite time singularity from a smooth initial 
condition \cite{Clay}. 
Global existence and regularity of the Navier-Stokes equations have 
been known in two space dimensions for a long time \cite{Lady70}. One 
of the main difficulties in obtaining the global regularity of the 
3D Navier-Stokes equations is mainly due to the presence of the vortex 
stretching, which is absent for the 2D problem. Under suitable 
smallness assumption on the initial condition, the local-in-time
existence and regularity results have been obtained for some 
time \cite{Lady70,Temam01,BM02}. But these results do not give
any hint on the question of global existence and regularity for
the 3D Navier-Stokes.

In this paper, we prove the nonexistence of local self-similar
singular solutions of the 3D Navier-Stokes equations. The local 
self-similar solutions we consider are very different from
the global self-similar solutions considered by Leray \cite{Leray34}.
The self-similar scaling is only valid in an inner core region which 
shrinks to a point dynamically as the time, $t$, approaches the 
singularity time, $T$. The solution outside the inner core region 
is assumed to be regular and does not satisfy self-similar scaling.
This type of local self-similar solution is developed dynamically,
and has been observed in some numerical studies. Under the assumption 
that the local self-similar velocity profile converges to a limiting 
profile as $t \rightarrow T$ in $L^p$ for some $p \in (3,\infty )$, 
we prove that such local self-similar blow-up is not possible. We 
remark that the nonexistence of global self-similar solutions has 
been proved by Necas, Ruzicka and Sverak in \cite{NRS96}. The result 
of \cite{NRS96} was further improved by Tsai in \cite{Tsai98}.

We prove our main result by using a ``Dynamic Singularity Rescaling'' 
technique. This technique is simple but effective. Below we give a 
brief description of this technique. Assume that the solution 
of the 3D Navier-Stokes develops a local self-similar singularity 
at $x=0$ at time $T$ for the first time. A typical local self-similar 
singular solution is of the form
\begin{equation}
    u ( x, t ) = \frac{1}{\sqrt{T - t}} U ( y, t), \quad
    p ( x, t ) = \frac{1}{T - t} P ( y, t), \quad 
    y = \frac{x}{\sqrt{T - t}},
\label{eqn-res0}
\end{equation}
for $0 \leqslant t < T$. We assume that $u$ is smooth outside an inner 
core region $\{x, |x | \leqslant (T-t)^\alpha \}$ for some $\alpha >0$ 
small. In particular, $u(x,t)$ and $p(x,t)$ are bounded for any 
fixed $|x| >0$ as $t \rightarrow T$. Using this condition, we can 
easily show that
\begin{equation}
| U(y, t) | \leqslant C (T)/|y|, \quad |P(y,t)| \leqslant C(T)/|y|^2 , \quad
\mbox{for } \;\; |y| \gg 1, \;\; t \leqslant T.
\label{eqn-decay}
\end{equation}
Thus, it is reasonable to assume that $U \in L^p$ for some 
$p \in (3,\infty )$.  But the $L^p$ norm of $U$ may be unbounded 
for $0 < p \leqslant 3$.

We assume that there exists a limiting profile $\overline{U}(y) \in L^p$ 
as $t \rightarrow T$ 
\begin{equation}
\lim_{t \rightarrow T} \| U(t) -  \Ub \|_{L^p} = 0 , 
\label{eqn-limit}
\end{equation}
for some $p$ satisfying $3 < p < \infty$.

Next, we introduce the following rescaling in time:
\begin{equation}
 \tau = \frac{1}{2} \log \frac{T}{T - t},
\label{eqn-restime}
\end{equation}
for $0 \leqslant t < T$. Note that by this time rescaling, we have 
transformed the original Navier-Stokes equations from $[0,T)$ in $t$ 
to $[0,\infty)$ in the new time variable $\tau$. Since 
$u$ is smooth for $0 < t < T$, $U$ is smooth for
$0 < \tau < \infty$. It is easy to derive the equivalent evolution 
equations for the rescaled velocity:
\begin{equation}
    U_\tau + U + (y \cdot \nabla) U + 2(U \cdot
    \nabla ) U = - 2 \nabla P + 2 \Delta U,
\label{eqn-res1}
\end{equation}
with initial condition $U|_{\tau = 0} = \sqrt{T} u_0 (y/\sqrt{T})$, where 
$U$ satisfies $\nabla \cdot U =0$ for all times.
The problem on the possible local self-similar 
blowup of the Navier-Stokes equations is now 
converted to the problem on the large time behavior of the rescaled 
equations (\ref{eqn-res1}). By assumption (\ref{eqn-limit}),
we know that $U \rightarrow \overline{U}$ as $\tau \rightarrow \infty$ 
in $L^p$. We will prove that the limiting velocity profile
actually satisfies the steady state equation of (\ref{eqn-res1}):
\begin{equation}
    \Ub + (y \cdot \nabla) \Ub + 2(\Ub \cdot
    \nabla ) \Ub = - 2 \nabla \Pb + 2 \Delta \Ub,
\label{eqn-res2}
\end{equation}
for some $\Pb$. Now it follows from
the result of \cite{Tsai98} that $\Ub \equiv 0$, which implies that
$\lim_{\tau \rightarrow \infty} \| U (\tau) \|_{L^p} = 0$ for some 
$p \in (3, \infty )$.

The fact that $\lim_{\tau \rightarrow \infty} \| U (\tau) \|_{L^p} = 0$
is significant because it shows that the rescaled velocity field becomes 
small dynamically as $\tau \rightarrow \infty$. It is easy to show that 
if the the solution $U$ is small in the $L^p$ norm at some time, $\tau_m$, 
the solution must decay exponentially in $\tau$ for $\tau \geqslant \tau_m$. 
The exponential decay in $U$ gives a sharp dynamic growth estimate in terms 
of the original velocity field. In fact, it exactly cancels the dynamic 
singular rescaling factor, $(\sqrt{T-t})^{-1}$, in the front of $U$. This 
gives us a uniform bound on the growth of $L^p$ for $0 < t < T$ with 
$p \in (3,\infty )$, and consequently it rules out the possibility of a 
finite time blowup at $T$ \cite{Prodi59,Serrin63,Lady67}. 

The nonexistence of local self-similar blowup of the 3D Navier-Stokes
equations has some interesting implication. First, the assumption
on the existence of a limiting self-similar profile, $\overline{U}$, 
can be verified numerically if a local self-similar blow is observed in 
a computation. Secondly, this result is related to a recent existence
result by one of the authors \cite{Hou06} for the axisymmetric 3D 
Navier-Stokes equations with swirl. Let $v^r$ denote the radial component
of the velocity field and $r=\sqrt{x^2 + y^2}$. 
The result in \cite{Hou06} shows that if 
$ \lim_{r \rightarrow 0} | r v^{r} | = 0$ holds uniformly
for $ 0 \leqslant t \leqslant T$, then the solution of the Navier-Stokes 
equations is regular for $t \leqslant T$. By the well-known 
Caffarelli-Kohn-Nirenberg result \cite{CKN82}, if the axisymmetric
3D Navier-Stokes equations develop a finite time singularity,
the singularity must lie in the $z$ axis. One of the most likely scenarios
that would violate the condition, 
$ \lim_{r \rightarrow 0} | r v^{r} | = 0$, is the local self-similar blowup of 
the Navier-Stokes equations. The result presented in this paper would
rule out such a possibility. For more discussions regarding other aspects 
of the Navier-Stokes equations, we refer the reader to 
\cite{Lady70,CM93,Temam01,BM02}.

The rest of the paper is organized as follows. In Section 2, we
state our main theorem and present its proof. The proof is 
divided into four subsections. 
A couple of technical results are deferred to the appendices.

\section{The main result and regularity analysis.}

\noindent
\begin{thm}\label{theorem1}

{\it Let $u_0 \in L^2(\mathbb{R}^3)\cap L^p(\mathbb{R}^3)$
for some $p \in (3, \infty )$ and $T$ be the first local self-similar
singularity time. Assume that $U(y,t)$ defined by (\ref{eqn-res0}) 
converges to $\Ub$ in $L^p$ as $t \rightarrow T$.
Then we must have $T = +\infty$, i.e. there is no finite time local
self-similar blowup for the 3D Navier-Stokes equations.}

\end{thm}

Before we prove our main theorem, we state the following well-known
$\left( L^{\tilde{q}}, L^{\tilde{p}} \right)$-estimates for the heat 
kernel in $\mathbb{R}^3$, 
$e^{- t \Delta}$, where $\Delta$ is the Laplacian operator.
  \begin{equation}
    \| e^{- t \Delta} w \|_{L^{\tilde{q}}} \leqslant c_0\; t^{-
    \left( \frac{3}{\tilde{p}} - \frac{3}{\tilde{q}} \right) / 2} \| w \|_{L^{
     \tilde{p}}},
\label{heat1}
  \end{equation}
\vspace{0.02in}
  \begin{equation}
    \| \nabla e^{- t \Delta} w \|_{L^{\tilde{q}}} \leqslant c_0\; 
    t^{- \left( 1 + \frac{3}{\tilde{p}} - \frac{3}{\tilde{q}} \right) / 2} 
   \| w \|_{L^{\tilde{p}}},
   \label{heat2}
  \end{equation}
for $1 < \tilde{p} \leqslant \tilde{q} < \infty$, $c_0$ depends on $\tilde{p}$ 
and $\tilde{q}$ only. In our analysis, we take $\tilde{q}=p$ and $\tilde{p}=p/2$. 
For this particular choice of $\tilde{p}$ and $\tilde{q}$, we can choose a constant, 
$c_0$, such that the above two inequalities hold. Throughout the paper, we will 
use $c_0$ and $c_1$ to denote various constants that do not depend on the 
individual functions, and use $C_j$ ($j=1,2$) to denote various constants 
that depend on the initial condition, $u_0$. We also define
\begin{equation}
\gamma = 3/p .
\label{gamma-def}
\end{equation}
Since $ 3 < p < \infty $, we have $ 0< \gamma < 1$. 

\vspace{0.2in}
\noindent
{\bf Proof of Theorem 1.}

We will prove the theorem by contradiction. Suppose that
$T < + \infty$. This means that the solution $u$ to problem
(\ref{nse}) develops a singularity at $t=T$ for the first time,
but $u$ is the unique smooth solution of (\ref{nse}) for $0 < t < T$
and is bounded in $L^p$. 

We will divide the proof into four steps, which are given
in the following four subsections.

%\vspace{0.1in}
%\noindent
\subsection{Dynamic singularity rescaling and a priori estimates.}
%\vspace{0.1in}

We make the following dynamic singularity rescaling of the 3D Navier-Stoke
equations:
\begin{equation}
  \left\{ \begin{array}{l}
    \tau = \dfrac{1}{2} \log \dfrac{T}{T - t}, \quad y = \dfrac{x}{\sqrt{T-t}} \\[2mm]
    u ( x, t ) = \dfrac{1}{\sqrt{T - t}} U ( y, \tau ),\\[2mm]
    p ( x, t ) = \dfrac{1}{T - t} P ( y, \tau ), \quad \mbox{for} \;\; 
0 \leqslant t < T .
  \end{array} \label{rescale} \right.
\end{equation}
Note that with this dynamic singularity rescaling, we transform 
the time interval from $[0,T)$ in the original time variable
$t$ to $[0,\infty)$ in the rescaled time variable $\tau$.  
It is easy to derive an evolution equation for the rescale velocity field:
\begin{equation}
  \left\{ \begin{array}{l}
    U_\tau + U + (y \cdot \nabla ) U + 2 (U \cdot
    \nabla ) U = - 2 \nabla P + 2 \Delta U, \\[2mm]
    \nabla \cdot U = 0 ,\\[2mm]
    U|_{\tau = 0} = \sqrt{T} u_0 \left( x \right) .
  \end{array} \label{U} \right.
\end{equation}

Note that since $u(x,t)$ is the unique smooth solution of the Navier-Stokes 
equations (\ref{nse}) for $0 < t < T$, $U(x,\tau)$ is the unique
smooth solution of the rescaled Navier-Stokes equations 
(\ref{U}) for $0 < \tau < \infty$.

Let $\phi (y) = (\phi_1,\phi_2,\phi_3)$ be a smooth, compactly 
supported, divergence free vector field in $\mathbb{R}^3$ and 
$\psi (\tau)$ be a smooth, compactly supported test function in
$(0,1)$ satisfying $\int_0^1 \psi (\tau) d \tau = 1$.
Multiplying (\ref{U}) by $\psi(\tau-n) \phi(y)$ and integrating
over $\mathbb{R}^3 \times [n,n+1]$ for some $n>0$, we obtain  
after integration by parts
\begin{eqnarray}
  && \int_n^{n+1} \int_{\mathbb{R}^3} \left ( -\psi_\tau \phi \cdot U  
 + \psi \phi \cdot U - \psi \nabla \cdot (\phi \otimes y ) \cdot U 
 - 2\psi \nabla \phi \cdot (U \otimes U) \right) dy d\tau \nonumber\\
  && = 2 \int_n^{n+1} \int_{\mathbb{R}^3}\psi \Delta \phi \cdot U dy d \tau ,
  \label{U-weak0}
\end{eqnarray}
where $\psi$ is evaluated at $\tau-n$.

By assumption (\ref{eqn-limit}), we have 
\begin{equation}
\lim_{\tau \rightarrow \infty} \| U(\tau ) - \Ub \|_{L^p} = 0,
\label{eqn-limit1}
\end{equation}
for some $p>3$. Thus $\| U (\tau) \|_{L^p}$ is bounded for $\tau$
sufficiently large, and $\| \Ub \|_{L^p}$ is also bounded. 
Let $U(\tau) = \Ub + R (\tau)$. By (\ref{eqn-limit1}), we have 
$\lim_{\tau \rightarrow \infty} \| R (\tau) \|_{L^p} = 0$.
Substituting $U(\tau) = \Ub + R (\tau)$ into (\ref{U-weak0}) and
letting $n \rightarrow \infty$, we will show that all the terms
involving $R$ will go to zero as $n\rightarrow \infty$. It is
sufficient to prove this for the nonlinear term:
\[
\int_n^{n+1} \int_{\mathbb{R}^3} \psi \nabla \phi \cdot (R \otimes R)
dy d \tau .
\]
Let $q = p/(p-2) > 1$. Then we have $2/p+1/q = 1$. Using the 
H\"older inequality, we obtain
\begin{eqnarray*}
\vert \int_n^{n+1} \int_{\mathbb{R}^3} \psi \nabla \phi \cdot (R \otimes R)
dy d \tau \vert &\leqslant & 
C \sup_{n \leqslant \tau \leqslant n+1} \int_{\mathbb{R}^3}
|\nabla \phi | |R |^2  d y \\
& \leqslant & C \| \nabla \phi \|_{L^q} \sup_{n \leqslant \tau \leqslant n+1}
\| R(\tau)\|_{L^p}^2 \rightarrow 0, \quad \mbox{as}\;\; n\rightarrow \infty.
\end{eqnarray*}
Other terms can be proved similarly. Therefore, by letting 
$n\rightarrow \infty$, we get 

\begin{eqnarray}
  && -\left (\int_0^{1} \psi_\tau (\tau) d \tau \right ) 
   \int_{\mathbb{R}^3} \phi (y) \Ub (y) dy \nonumber \\
 && + \left (\int_0^{1} \psi (\tau) d \tau \right ) 
    \left ( \int_{\mathbb{R}^3} \left ( \phi \cdot \Ub - 
     \nabla \cdot (\phi \otimes y ) \cdot \Ub
 - 2\nabla \phi \cdot (\Ub \otimes \Ub ) \right) dy \right ) \nonumber\\
  && = 2 \left (\int_0^{1} \psi (\tau ) d \tau \right ) 
   \left ( \int_{\mathbb{R}^3} \Delta \phi \cdot \Ub dy \right ).
  \label{U-weak}
\end{eqnarray}
Since $\psi $ has compact support in $[0,1]$, we conclude that
\[
\int_0^{1} \psi_\tau (\tau) d \tau  = 0.
\]
Moreover, we have $\int_0^{1} \psi (\tau ) d \tau = 1$ by assumption
on $\psi$. Thus, we obtain

\begin{eqnarray}
\int_{\mathbb{R}^3} \left ( \phi \cdot \Ub -
     \nabla \cdot (\phi \otimes y ) \cdot \Ub
 - 2\nabla \phi \cdot (\Ub \otimes \Ub ) - 2\Delta \phi \cdot \Ub \right ) dy = 0.
  \label{U-weak1}
\end{eqnarray}
Thus, $\Ub$ is a weak solution of the steady state rescaled 
Navier-Stokes equations:
\begin{eqnarray}
\Ub + (y \cdot \nabla) \Ub + 2 (\Ub \cdot \nabla )\Ub =  - 2 \nabla \Pb + 
2\Delta \Ub , 
\label{U-steady}
\end{eqnarray}
with $\nabla \cdot \Ub = 0$. Let $R_j$ be a Riesz operator with Fourier 
symbol $\xi_j/|\xi |$. One can easily modify the proof of Lemma 3.1 of 
\cite{NRS96} to show that $\Pb = R_jR_k (\Ub_j \Ub_k)$. 

Since $\Ub \in L^p $ for some $p \in (3,\infty)$, 
we can apply Theorem 1 of \cite{Tsai98} 
to conclude that $\Ub \equiv 0$. As a result, we obtain the following
{\it a priori} decay estimate for $\| U (\tau) \|_{L^p}$.

\begin{lemma}
\label{lemma1}  
The solution $U(x,\tau)$ of the rescaled Navier-Stokes equations 
(\ref{U}) satisfies the following decay estimate:
  \begin{equation}
      \lim_{\tau \rightarrow \infty} \| U (\tau ) \|_{L^p}  = 0.
     \label{Ulp-decay}
  \end{equation}
\end{lemma}

For the purpose of our later analysis, we will choose a $\tau_m$ 
large enough to satisfy the following inequality:
\begin{equation}
2 c_0^2\; c_1 \| U (\tau_m)  \|_{L^p} \leqslant \frac{1}{6},
\label{tau-m}
\end{equation}
where the constant $c_1$ is defined by
\begin{equation}
c_1 = (\frac{2}{1-\gamma} + \frac{1}{2}) \left(1-e^{-2}\right)^{- \frac{1+\gamma}{2}}.
\label{C1-def}
\end{equation}
The reason for such a choice of $\tau_m$ will become clear later in our analysis.

\subsection{Dynamic decay estimates for the rescaled equations.}

In this subsection, we will perform estimates for the
rescaled Navier-Stokes equations starting from $\tau = \tau_m$ 
with the initial value, $U(x,\tau_m)$:
\begin{equation}
  \left\{ \begin{array}{l}
    V_\tau + V + (y \cdot \nabla ) V + 2  (V \cdot \nabla) V = - 2  \nabla P + 2 \Delta V\\[2mm]
    \nabla \cdot V = 0\\[2mm]
    V|_{\tau = 0}(x) \equiv V_0 (x) = U(x,\tau_m),
  \end{array} \label{NSE-rescaled} \right.
\end{equation}
where $\tau_m$ is defined by (\ref{tau-m})-(\ref{C1-def}).
Since $U(x,\tau)$ is the unique smooth solution of the rescaled Navier-Stokes 
equations (\ref{U}) for $0 < \tau < \infty$, the function $V(x,\tau)$
defined by
\begin{equation}
V(x,\tau) = U(x,\tau+\tau_m), \quad \mbox{for} \;\; \tau \geqslant 0,
\label{V-U}
\end{equation}
is the unique smooth solution of (\ref{NSE-rescaled}).

Next, we perform estimates for the linearized operator in (\ref{NSE-rescaled})
\begin{equation}
  \frac{\partial V}{\partial \tau} + V + ( y \cdot \nabla_y ) V - 2 \Delta_y V = 0,
  \label{lin-op0}
\end{equation}
with initial value $V|_{\tau = 0} = V_0 $. 

Let $y = e^\tau {\tilde y}$ and $\tilde{V}(\tilde y,\tau) \equiv V(y,\tau)$. 
Then we have
\begin{equation}
  \frac{\partial \tilde{V}}{\partial \tau} + \tilde{V} - 2 e^{-2\tau} 
\Delta_{\tilde y} \tilde{V} = 0,
  \label{lin-op}
\end{equation}
with initial value $\tilde{V}|_{\tau = 0} = V_0 $.

Taking the Fourier transform of (\ref{lin-op}), we get
\begin{equation}
  \frac{\partial \widehat{\tilde V}}{\partial \tau} + \widehat{\tilde V} + 2 e^{- 2 \tau} | \xi
  |^2 \widehat{\tilde V} = 0,
\label{eqn:fourier}
\end{equation}
where the Fourier transformation is defined as $\widehat{f} ( \xi ) \equiv
\int f ( x ) e^{- 2 \pi i x \cdot \xi} d x$. Equation (\ref{eqn:fourier}) can be 
written as
\begin{equation}
  \frac{\partial}{\partial \tau} \left( e^{\tau + 2 | \xi |^2 \int_0^{\tau}
  e^{- 2 s} d s} \widehat{\tilde V} ( \tau ) \right) = 0.
\end{equation}
Integrating from $0$ to $\tau$, we get
\begin{equation}
    \widehat{\tilde V} ( \tau ) 
     =  e^{- \tau - | \xi |^2 \left( 1 - e^{- 2 \tau} \right)} \widehat{V_0} .
\end{equation}
Using the explicit formula of the Fourier transform of a Gaussian in three
space dimensions 
(see, e.g. \cite{Stein70}) 
\begin{equation}
  \widehat{e^{- \pi \alpha^2 |x|^2}} = \frac{1}{\alpha^3} e^{- \pi | \xi |^2 /
  \alpha^2},
\end{equation}
with $\alpha^2 = \left( \frac{\pi}{ 1 - e^{- 2 \tau} } \right)$,
we obtain 
\begin{equation}
 {\cal F}^{-1} \left (e^{- | \xi |^2 ( 1 - e^{- 2 \tau} ) }\right ) 
= \left(
  \frac{\pi}{ 1 - e^{- 2 \tau} } \right)^{\frac{3}{2}} e^{-
  \pi^2 |x|^2 /  ( 1 - e^{- 2 \tau} )},
\end{equation}
where ${\cal F}^{-1}{f} ( x ) \equiv \int f ( \xi ) e^{ 2 \pi i x \cdot \xi} d \xi$
is the inverse Fourier transformation. Therefore, we have
\begin{equation}
  {\tilde V} ( \tilde y , \tau ) = e^{-\tau} \left( \frac{\pi}{ 1 - e^{- 2 \tau} }
  \right)^{\frac{3}{2}} \int V_0 \left( \tilde{x} \right)
  \left ( e^{- \pi^2 | \tilde{y} -
  \tilde{x} |^2 / ( 1 - e^{- 2 \tau} )}\right ) d \tilde{x} . \label{anay-sol}
\end{equation}
Denote by $e^{- \tau A}$ the solution operator of the linearized equations
(\ref{lin-op}). Define
\begin{equation}
t_0 (\tau ) = ( 1 - e^{- 2 \tau} ),
\label{tau-0}
\end{equation}
and denote $\Delta$ as the Laplacian operator, then we have
\begin{equation}
  e^{- \tau A} V_0 = \tilde {V} (\tilde y,\tau ) = e^{- \tau} 
\left (e^{- t_0 ( \tau ) \Delta} V_0 \right ).
\label{A-operator}
\end{equation}

Define the following bilinear operator:
\begin{equation}
  F ( U, V ) = 2  \left( 1 - \nabla \left( - \Delta
  \right)^{- 1} \nabla \cdot \right) \nabla \cdot \left( U \otimes V \right).
  \label{fuv}
\end{equation}
In particular, if we set $V=U$, we have
\begin{equation}
  \begin{array}{lll}
  F ( U, U ) & = & 2 \left( \nabla \cdot \left( U     
    \otimes U \right) - \nabla \left( - \Delta \right)^{- 1} \nabla
    \cdot \nabla \cdot \left( U \otimes U \right) \right)\\[2mm]
    & = & 2 \left( U \cdot \nabla U + \nabla P \right).
  \end{array} \label{fu}
\end{equation}

The rescaled 3D Navier-Stokes equations
(\ref{NSE-rescaled}) can be converted into the following integral equation: 
\begin{equation}
V(\tau) = e^{-\tau A} V_0 - \int_0^{\tau} e^{- \left( \tau - s \right) A} 
F ( U, U ) ( s ) d s .
\label{NSE-integral}
\end{equation}

To solve the integral equation (\ref{NSE-integral}), we 
construct a successive approximation, $V^{(n)}$, using the
following iterative scheme (see \cite{Kato84}): 
$V^{(0)} =  e^{-\tau A} V_0$,
\begin{equation}
  V^{(n + 1)} = V^{(0)} - G ( V^{(n)},V^{(n)}), \quad
 n \geqslant 0,
\label{it-scheme}
\end{equation}
where the bilinear operator $G(U,V)$ is defined as follows:
\begin{equation}
  G ( U, V ) = \int_0^{\tau} e^{- \left( \tau - s \right) A} F ( U, V ) ( s )
  d s .
\end{equation}

To establish the convergence of the approximate solution sequence,
$ V^{(n)}$, we need to use the following lemma, which follows
from (\ref{A-operator}) and the well-known 
$\left( L^q, L^p \right)$-estimates (\ref{heat1})-(\ref{heat2}) 
for the heat kernel. 

\begin{lemma}

Let $V \in L^{\tilde{p}} $ for $1 <\tilde{p} \leqslant \tilde{q} < \infty$. We have
  \begin{equation}
    \| e^{- \tau A} V \|_{L^{\tilde{q}}} \leqslant c_0\; e^{- (1-3/\tilde{q})\tau} t_0 ( \tau )^{-
    \left( \frac{3}{\tilde{p}} - \frac{3}{\tilde{q}} \right) / 2} 
   \| V \|_{L^{\tilde{p}}} ,\label{l-pq}
  \end{equation}
  \begin{equation}
    \| \nabla e^{- \tau A} V \|_{L^{\tilde{q}}} \leqslant c_0\;e^{- (2-3/\tilde{q})\tau} t_0 ( \tau
    )^{- \left( 1 + \frac{3}{\tilde{p}} - \frac{3}{\tilde{q}} \right) / 2} 
    \| V \|_{L^{\tilde{p}}} .
    \label{dl-pq}
  \end{equation}
\end{lemma}
The lemma can be proved easily by noting that the heat kernel actually 
acts on the variable $\tilde y$ through the function ${\tilde V}(\tilde y ,\tau)$ and 
$\tilde y = e^{-\tau }y$. Thus we lose a factor $e^{3\tau /{\tilde q}}$ when we
estimate the $L^{\tilde q}$ norm by changing variables from $y$ to $\tilde y$, but we
gain a factor of $e^{-\tau}$ when we differentiate with respect to $y$. 

Applying (\ref{l-pq}) with $\tilde{p}=\tilde{q}=p$, we obtain 
\begin{equation}
\| V^{(0)} \|_{L^p} ( \tau ) = \| e^{-\tau A} V_0 \|_{L^p} ( \tau )
\leqslant c_0 e^{- (1-\gamma )\tau} \| V_0 \|_{L^p},
\label{est-V0}
\end{equation}
where $\gamma = 3/p$. To estimate $ \| G ( U, V ) \|_{L^p}$, we use (\ref{dl-pq})
with $\tilde{q}=p$ and $\tilde{p}=p/2$: 
\begin{equation}
  \| G ( U, V ) \|_{L^p} ( \tau ) \leqslant 2 c_0 
  \int_0^{\tau} e^{- (2-\gamma) ( \tau - s )} t_0 ( \tau - s )^{- \frac{1+\gamma}{2}}
  \| U \|_{L^p} (s)\| V \|_{L^p}(s) d s \label{guv},
\end{equation}
where we have used 
the H\"older inequality
$\| U \otimes V \|_{L^{p/2}} \leqslant \| U \|_{L^p} \| V \|_{L^p}$
and
the fact that $\left( - \Delta \right)^{- 1} \nabla \cdot
\nabla \cdot$ is a Rietz operator of degree zero, which is a bounded
operator from $L^p$ to $L^p$. 
In particular, we obtain by setting $V=U$ that
\begin{equation}
  \| G ( U,U ) \|_{L^p} ( \tau ) \leqslant 2 c_0 
  \int_0^{\tau} e^{- (2-\gamma)( \tau - s )} t_0 ( \tau - s )^{- \frac{1+\gamma}{2}}
   \| U \|_{L^p}^2 (s) d s \;. \label{gu}
\end{equation}
Now, applying (\ref{est-V0}) and (\ref{gu}) to the iterative scheme 
(\ref{it-scheme}), we get
\begin{equation}
  \| V^{(n + 1)} \|_{L^p} ( \tau ) \leqslant c_0 e^{-(1-\gamma)\tau} \| V_0 
  \|_{L^p} + 2c_0 \int_0^{\tau}  e^{- (2-\gamma)(\tau-s)} t_0 \left( \tau -
  s \right)^{- \frac{1+\gamma}{2}} \| V^{(n)} \|_{L^p}^2 ( s ) d s.
  \label{it-norm}
\end{equation}
Define
\begin{equation}
  K_n = \sup_{0 \leqslant \tau < \infty} \| e^{(1-\gamma)\tau} V^{(n)} ( \tau )
  \|_{L^p} . \label{kn}
\end{equation}
Multiplying (\ref{it-norm}) by $e^{(1-\gamma)\tau}$ on both sides and
using (\ref{kn}), we obtain
\begin{equation}
    e^{(1-\gamma)\tau} \| V^{(n + 1)} \|_{L^p} ( \tau ) 
     \leqslant c_0 \| V_0 \|_{L^p} + 2 c_0 e^{-\tau}
    K_n^2 \int_0^{\tau} e^{\gamma s}
     \left( 1 - e^{- 2 ( \tau - s )} \right)^{- \frac{1+\gamma}{2}}
    ds .
   \label{it-norm1}
\end{equation}
In Appendix II,
we will prove that 
\begin{equation}
\label{C1-integral}
e^{-\tau}\int_0^{\tau} e^{\gamma s} \left( 1 - e^{- 2 ( \tau - s )} 
\right)^{- \frac{1+\gamma}{2}} d s \leqslant c_1, \quad \mbox{for}
\;\;\mbox{all}\;\; \tau \geqslant 0,
\end{equation}
where $c_1$ is defined in (\ref{C1-def}). 
Now, take the supremum of the both sides of (\ref{it-norm1}) for all 
$\tau \geqslant 0$, we obtain the following recurrence inequalities:
\begin{equation}
  K_{n + 1} \leqslant K_0 + M K_n^2 , \quad \mbox{for} \;\; 
n\geqslant 0,\label{ind-ine}
\end{equation}
with $K_n |_{n=0} = K_0$,
where 
\begin{equation}
M = 2 c_0 c_1, \quad
K_0 = c_0 \| V_0 \|_{L^p}.
\label{M-def}
\end{equation} 
We will prove the following lemma in Appendix I.

\begin{lemma}
  \label{lemma-ind-ine} Let $K_0$ and $M$ be two positive
constants satisfying
  \begin{equation}
    K_0 M \leqslant \frac{1}{6}, \label{k0c1}
  \end{equation}
  then there exists a positive constant $K_{\max}$, such that
  \begin{equation}
    K_n \leqslant K_{\max}, \hspace{2em} \tmop{for} \;\; \tmop{all} \;\; 
n \geqslant 1
    \label{kn-bound}
  \end{equation}
  holds for the recurrence sequence $K_n$ satisfying (\ref{ind-ine}). 
  Moreover the upper bound $K_{\max}$ satisfies
  \begin{equation}
    2 M K_{\max} \leqslant \frac{1}{2} .
\label{2mkn}
  \end{equation}
\end{lemma}

Recall that in (\ref{tau-m}), we have chosen $\tau_m$ such that
\begin{equation}
  2 c_0^2 c_1 \| U ( \tau_m ) \|_{L^p} \leqslant 1/6 \;.
\end{equation}
Therefore, we have
\begin{equation}
K_0 M \leqslant
2 c_0^2 \;c_1 \;\| U ( \tau_m ) \|_{L^p}
\leqslant \frac{1}{6} . 
\end{equation}
Thus, for our choice of $\tau_m$ defined in (\ref{tau-m}), the recurrence
sequence $K_n$ has an upper bound $K_{\max}$ for all $n$. That is
\begin{equation}
\| V^{(n)} \|_{L^p} (\tau) \leqslant K_{\max} e^{-(1-\gamma)\tau} , 
\quad \mbox{for} \quad n \geq 1.
\label{Vn-bound}
\end{equation}

\subsection{Convergence of the approximate solution sequence.}

In this subsection, we will establish the convergence of the 
approximate solution sequence, $\{V^{(n)}\}$, and study the 
property of its limiting solution. We will first 
show that the approximate solution sequence
$\{V^{(n)}\}$ is a Cauchy sequence in $L^p$. By
subtracting (\ref{it-scheme}) with index $n$ from that with index
$n - 1$, we obtain 
\begin{eqnarray*}
    && \| V^{(n + 1)} - V^{(n)} \|_{L^p} \\[2mm]
   && =  \| G ( V^{(n)},V^{(n)} ) - G (V^{(n-1)},V^{(n-1)}) \|_{L^p} \\[2mm]
   && = \| G ( V^{(n)}, V^{(n)} - V^{(n-1)} ) + G( V^{(n)} - V^{(n-1)}, V^{(n-1)})
    \|_{L^p}.
\end{eqnarray*}
Using (\ref{guv}), (\ref{Vn-bound}) and (\ref{C1-integral}), we obtain
\begin{eqnarray*}
&& e^{(1-\gamma)\tau} \| V^{(n + 1)} - V^{(n)} \|_{L^p} ( \tau )\\[2mm]
 && \leqslant 2 c_0 e^{(1-\gamma)\tau}
  \int_0^{\tau} e^{- (2-\gamma)(\tau-s)} t_0 \left( \tau - s
    \right)^{- \frac{1+\gamma}{2}} \left( \| V^{(n)} \|_{L^p} 
    + \| V^{(n-1)} \|_{L^p} \right) \| V^{(n)} - V^{(n-1)}
    \|_{L^p} ( s ) d s\\[2mm]
   &&\leqslant 4  c_0 K_{\max} e^{-\tau} \int_0^{\tau} e^{\gamma s}t_0 \left( \tau - s
    \right)^{- \frac{1+\gamma}{2}} d s \left ( \sup_{0 \leqslant s <
    \infty} e^{(1-\gamma)s}\| V^{(n)} - V^{(n - 1)} \|_{L^p} (s)\right )\\[2mm]
    &&\leqslant 4  c_0 c_1 K_{\max} \sup_{0 \leqslant s < \infty} 
    e^{(1-\gamma)s}\| V^{(n)} -
    V^{(n - 1)} \|_{L^p} (s) \\[2mm]
    &&\leqslant \frac{1}{2} \sup_{0 \leqslant s < \infty}
    e^{(1-\gamma)s} \| V^{(n)} - V^{(n - 1)} \|_{L^p} (s),
\end{eqnarray*}
where we have used (\ref{M-def}) and (\ref{2mkn}) in deriving the last 
inequality. Taking the supremum on the left hand side would yield
\begin{equation}
  \sup_{0 \leqslant \tau < \infty} e^{(1-\gamma)\tau} \| V^{(n + 1)} -
  V^{(n)} \|_{L^p} (\tau) \leqslant \frac{1}{2} \sup_{0 \leqslant \tau
  < \infty}  e^{(1-\gamma)\tau} \| V^{(n)} - V^{(n - 1)} \|_{L^p} (\tau),
\end{equation}
which implies
\begin{equation}
\sup_{0 \leqslant \tau < \infty} e^{(1-\gamma)\tau} \| V^{(n + m)} -
V^{(n)} \|_{L^p} (\tau) \leq C_1 \left (\frac{1}{2}\right )^n,  
\quad \mbox{for} \;\; \mbox{any} \; n, m \geq 1,
\end{equation}
where $C_1$ depends on $V^{(0)}$ only. Thus $\{V^{(n)} \}$ 
is a Cauchy sequence in $\tmop{BC} \left( \left[ 0, \infty \right)
; L^p \left( \mathbb{R}^3 \right) \right)$. Here 
 $\tmop{BC} \left( \left[ 0, \infty \right)
; L^p \left( \mathbb{R}^3 \right) \right)$ denotes the class
of bounded and continuous function from $[ 0, \infty )$
to $L^p (\mathbb{R}^3 )$. As a result, we
have proved that $V^{(n)} \left( \tau \right)$ converges
uniformly to a limiting function $\overline{V} \left( \tau \right)$ in 
$\tmop{BC} \left( \left[ 0, \infty \right)
; L^p\left( \mathbb{R}^3 \right) \right)$. Taking the
limit $n \rightarrow \infty$ in (\ref{Vn-bound}), we 
obtain 
\begin{equation}
\| \overline{V} \|_{L^p} (\tau) \leqslant K_{\max} e^{-(1-\gamma)\tau} .
\label{V-bound}
\end{equation}

Next, we will show that $\overline{V}$ is a solution of
the integral equation (\ref{NSE-integral}). To this end, we define 
$R^{(n)}(x,\tau) \equiv V^{(n)}(x,\tau) - \overline{V}(x,\tau)$.
We have just shown that
\begin{equation}
\sup_{0 \leqslant \tau < \infty} e^{(1-\gamma)\tau} \| R^{(n)} \|_{L^p} (\tau) =
\sup_{0 \leqslant \tau < \infty} e^{(1-\gamma)\tau} 
\| V^{(n)}- \overline{V} \|_{L^p} (\tau) \rightarrow 0,
\label{Rn-error}
\end{equation}
as $n \rightarrow \infty$. 
Now substituting $V^{(n)} = \overline{V} + R^{(n)}$ into the iterative scheme
(\ref{it-scheme}) and using the bilinearity of operator $G(U,V)$, we get 
\begin{equation}
\overline{V} - V^{(0)} + G(\overline{V},\overline{V}) = -(R^{(n+1)} + G(R^{(n)},\overline{V}) + 
G(\overline{V},R^{(n)})+ G(R^{(n)},R^{(n)})) .
\label{eqn-error}
\end{equation}
We will prove that the error terms on the right hand side of 
(\ref{eqn-error}) tend to zero uniformly for all $\tau \geqslant 0$.
It is obvious that $\|R^{(n+1)}\|_{L^p} \rightarrow 0 $ 
uniformly as $n \rightarrow \infty$ from (\ref{Rn-error}). 

To show that the error terms which are linear in $R^{(n)}$ tend to zero 
uniformly, we use (\ref{guv}) and the {\it a priori} bound on 
$\overline{V}$ given by (\ref{V-bound}). Specifically, we have
\begin{eqnarray*}
&& \| G(R^{(n)},\overline{V}) + G(\overline{V},R^{(n)}) \|_{L^p} ( \tau )\\[2mm]
 && \leqslant 4 c_0 \int_0^{\tau} e^{- (2-\gamma)(\tau-s)} t_0 \left( \tau - s
    \right)^{- \frac{1+\gamma}{2}} \| R^{(n)} \|_{L^p} ( s )
    \|\overline{V} \|_{L^p} ( s ) d s\\[2mm]
   &&\leqslant 4  c_0 K_{\max} e^{-(1-\gamma)\tau} e^{-\tau} \int_0^{\tau}
 e^{\gamma s} t_0 \left( \tau - s \right)^{- \frac{1+\gamma}{2}} d s 
\left ( \sup_{0 \leqslant s < \infty} e^{(1-\gamma)s}\| R^{(n)}\|_{L^p}(s)\right )\\[2mm]
    &&\leqslant 4  c_0 c_1 e^{-(1-\gamma)\tau} K_{\max}
    \sup_{0 \leqslant s < \infty} 
    e^{(1-\gamma)s} \| R^{(n)} \|_{L^p}(s)\\[2mm] 
    &&\leqslant 
    \sup_{0 \leqslant s < \infty}e^{(1-\gamma)s} \| R^{(n)} \|_{L^p}
    (s) \rightarrow 0 ,
\end{eqnarray*}
uniformly for all $\tau$ as $n \rightarrow \infty$, where we have
used $M=2 c_0 c_1 $ and (\ref{2mkn}).

To show that the nonlinear error term 
$G(R^{(n)},R^{(n)})$ also tends to zero uniformly, we
note that the {\it a priori} bounds on $V^{(n)}$ and $\overline{V}$ also provide
the following {\it a priori} bound for $R^{(n)}$:
\begin{equation}
\| R^{(n)} \|_{L^p} (\tau) \leqslant 2K_{\max} e^{-(1-\gamma)\tau} ,
\quad \mbox{for} \quad n \geq 1.
\label{Rn-bound}
\end{equation}
Using (\ref{Rn-bound}) and applying the same argument as above, we can prove that
\[
 \| G(R^{(n)},R^{(n)}) \|_{L^p} ( \tau )
\leqslant \sup_{0 \leqslant s < \infty} e^{(1-\gamma)s} \| R^{(n)} \|_{L^p} \rightarrow 0,
\]
uniformly for $0 \leqslant \tau < \infty$ as $n \rightarrow \infty$. 

Now, passing the limit $n\rightarrow \infty$ in the $L^p$ norm, we obtain
\begin{equation}
\overline{V} (\tau)= V^{(0)} - G(\overline{V},\overline{V}) , 
\quad \mbox{for} \;\; \mbox{all} \;\; \tau \geqslant 0,
\label{eqn-Vbar}
\end{equation}
which shows that $\overline{V}$ is a solution of the integral equation 
(\ref{NSE-integral}), satisfying the decay property (\ref{V-bound}).

\subsection{The non-blowup estimates in the original variables.}

In this subsection, we will complete the regularity analysis
in the original physical variable.
By the uniqueness of strong solutions in $L^p$ with $p>3$,
we have
\begin{equation}
\|\overline{V} \|_{L^p}(\tau) = \| U \|_{L^p} (\tau +\tau_m), 
\quad \mbox{for} \;\; 0 \leqslant \tau < \infty .
\label{VU-norm}
\end{equation}

%To prove (\ref{VU-norm}), we recall that $V(x,\tau ) \equiv U(x,\tau+\tau_m )$ 
%is the unique smooth solution of (\ref{NSE-rescaled}) and is bounded 
%in $L^p$. Since the original Navier-Stokes equations (\ref{nse}) for 
%$0 \leqslant t < T$ is equivalent to the rescaled Navier-Stokes equations 
%(\ref{NSE-rescaled}) for $0 \leqslant \tau < \infty$, the uniqueness of 
%strong solutions of (\ref{nse}) in $L^p$ with $p>3$ 
%\cite{Prodi59,Serrin63} implies 
%the uniqueness of the rescaled Navier-Stokes equations (\ref{NSE-rescaled}) 
%in $L^p$. Thus (\ref{VU-norm}) holds. Equality (\ref{VU-norm}) can also 
%be proved directly by a simple contraction mapping argument. For the sake of 
%completeness, we provide a simple proof of (\ref{VU-norm}) in Appendix III.
%
Now we can use the decay estimate for $\overline{V}$ in (\ref{V-bound})
to obtain a decay estimate for $U$, which in turn will rule out the
possibility of a finite time singularity for the 3D Navier-Stokes
equations.

Using (\ref{V-bound}) and (\ref{VU-norm}), we immediately obtain 
a decay estimate for $U$:
\begin{equation}
  \| U \|_{L^p} \left( \tau \right) \leqslant K_{\max}
  e^{-(1-\gamma)( \tau - \tau_m )} , \hspace{2em} \tmop{for} \;\; 
  \tau \geqslant \tau_m . \label{Ul4-decay1}
\end{equation}
This proves the following decay estimate for $U$.

\begin{lemma}
  The solution $U(x,\tau)$ of rescaled Navier-Stokes equations
(\ref{U}) with $\tau_m$ defined by (\ref{tau-m}) has
  a uniform decay rate in $\tau$ as follows:
  \begin{equation}
    \| U \|_{L^p} \left( \tau \right) \leqslant K_{\max} e^{- (1-\gamma)( \tau -
    \tau_m )}, \hspace{2em} \tmop{for}  \;\; \tau \geqslant \tau_m .
    \label{Ul4-decay2}
  \end{equation} \end{lemma}

Substituting the relation 
\begin{equation}
  u ( x, t ) = \frac{1}{\sqrt{T - t}} U ( y, \tau )
\end{equation}
into (\ref{Ul4-decay2}), we obtain for 
$t_m \leqslant t < T$ with $t_m = T \left( 1 - e^{- 2 \tau_m} \right)$,
\begin{equation}
  \begin{array}{lll}
    \| u \|_{L^p} ( t ) & = & \dfrac{(T - t)^{\gamma/2}}{(T - t)^{1/2}}
    \| U \|_{L^p} ( \tau )\\[3mm] 
    & \leqslant & \dfrac{K_{\max}}{(T - t)^{(1-\gamma)/2}} 
    e^{- (1-\gamma)( \tau - \tau_m )} = \dfrac{K_{\max} e^{(1-\gamma)\tau_m}}
    {(T - t)^{(1-\gamma)/2}} e^{- (1-\gamma)\tau}\\ [3mm]
    & = & \dfrac{K_{\max} e^{(1-\gamma)\tau_m}}{(T - t)^{(1-\gamma)/2}} 
     \left(\dfrac{T - t}{T}\right )^{(1-\gamma)/2}\\[3mm]
    & \leqslant & \dfrac{K_{\max} e^{(1-\gamma)\tau_m}}{T^{(1-\gamma)/2}}, \hspace{2em}
    \tmop{for} \quad t_m \leqslant t < T.
  \end{array}
\label{u-decay}
\end{equation}
Since $u_0 \in L^p$ for some $p \in (3,\infty)$, it is easy to show that
there is a local-in-time smooth solution whose $L^p$ norm is bounded
\cite{Kato84} (This can also be proved directly by using the same iterative 
scheme applied to the original Navier-Stokes equations for a short time). 
Moreover, since $T$ is the first singularity time, we conclude that
$u$ is smooth for $0 < t \leqslant t_m < T$ and has
a bounded $L^p$ norm for $ t \leqslant t_m$. Thus, $\| u \|_{L^p} (t)$ is 
uniformly bounded for $ 0 \leqslant t < T$.

Now, we can apply the so-called Ladyzhenskaya-Prodi-Serrin condition
(see \cite{Lady67}, \cite{Prodi59} and \cite{Serrin63}), which is also
known as the $L^{p,q}$ criteria. The so-called $L^{p,q}$ criteria state
that if a suitable weak solution of (\ref{nse}) satisfies
\begin{equation}
  u \in L^q \left( [ 0, T ) ; L^p \left( \mathbb{R}^3 \right) \right) 
\end{equation}
with
\begin{equation}
  \frac{3}{p} + \frac{2}{q} \leqslant 1, \hspace{2em} p \in [ 3, \infty ],
\end{equation}
then $u$ is a smooth solution of the 3D Navier-Stokes equation
up to $t = T$. In our case, we have obtained a uniform bound
in $L^p$ for $u$ with $p \in (3,\infty )$ for $0 \leqslant t < T$. Thus the 
$L^{p,q}$ criterion is satisfied with $q = \infty$. Therefore,
we conclude that $u$ is a smooth function in 
$\mathbb{R}^3 \times ( 0, T ]$.

This conclusion contradicts with our assumption that $u$ would cease
to be regular at time $T$ for the first time. This contradiction implies
that $u$ can not develop a local self-similar singularity in any finite time. 
This completes the proof of Theorem 1.

\vspace{0.2in}
\centerline{\bf \Large Appendix I.}

\vspace{0.2in}
In this appendix, we prove Lemma \ref{lemma-ind-ine}.

\noindent
{\bf Proof of Lemma \ref{lemma-ind-ine}}.
  It is sufficient to obtain an upper bound for the recurrence equalities
  \begin{equation}
    \widetilde{K}_{n + 1} = \tilde{K}_0 + M 
    \tilde{K}_n^2, \hspace{2em} \tilde{K}_0 = K_0 .
  \end{equation}
  It is easy to see that $K_n \leqslant \tilde{K}_n$, for all $n \geqslant 1$. To
  simplify the notation, we will drop the tilde in $\tilde{K}_n$ in the
  following. Define $l_n = K_{n + 1} - K_n$, then
  we have
  \begin{equation}
    l_{n} = M \left( K_{n - 1} + K_n \right) l_{n-1}.
   \label{ln-def}
  \end{equation}
  It is easy to see that $l_n >0$ for all $n \geqslant 0$ and
  $K_n$ is a monotonely increasing sequence. We claim that
  \begin{equation}
    M \left( K_{j - 1} + K_j \right) \leqslant \frac{1}{2}, \hspace{2em}
    \tmop{for} \;\; \tmop{all} \;\; j \geqslant 1. \label{claim}
  \end{equation}
We will prove (\ref{claim}) by an induction argument.
  \begin{enumerate}
    \item For $j = 1$, we have
    \begin{equation}
      M \left( K_0 + K_1 \right) = M \left( K_0 + K_0 + M K_0^2 \right)
      \leqslant \frac{1}{2}
    \end{equation}
    from the assumption $K_0 M \leqslant \frac{1}{6}$.
    
    \item Assume that (\ref{claim}) holds for all $j \leqslant n$, we will
     prove that it also hold for $j=n+1$. Let $\alpha = 1/2$. It 
     follows from (\ref{ln-def}) and the induction assumption that
    \begin{equation}
      l_j \leqslant \alpha l_{j - 1}, \quad \mbox{for}\;\;\mbox{all}\;\;
      1 \leqslant j \leqslant n,
    \end{equation}
     which implies that
    \begin{equation}
     l_j \leqslant \alpha^j l_0 .
     \label{lj-bound}
    \end{equation}
     Applying $K_{n+1} = K_n + l_n$ recursively and using (\ref{lj-bound}), 
     we obtain
    \begin{equation}
      \begin{array}{lll}
        K_{n + 1} & = & K_0 + \sum_{j=0}^n l_{j}\\[2mm]
        & \leqslant & K_0 + l_0 \sum_{j=0}^n \alpha^j\\[2mm]
        & = & K_0 + l_0 \dfrac{1 - \alpha^{n + 1}}{1 - \alpha}\\[2mm]
        & \leqslant & K_0 + 2 M K_0^2 \leqslant \frac{4}{3} K_0,
      \end{array}
    \end{equation}
    where we have used $M K_0 \leqslant 1/6$.
    Define $K_{\max} = \frac{4}{3} K_0$. Then we have 
\begin{equation}
2 M K_{\max} = \frac{8}{3} M K_0 \leqslant \frac{4}{9} < \frac{1}{2}.
\label{Kmax}
\end{equation}
    Thus, we obtain
    \begin{equation}
        M \left( K_n + K_{n + 1} \right)  \leqslant  2 M K_{\max}
        < \frac{1}{2} .
    \end{equation}
  \end{enumerate}
  This proves the claim (\ref{claim}) by induction, and we obtain
  \begin{equation}
    K_{n} \leqslant K_{\max},\quad \mbox{for}\;\;\mbox{all}\;\; n\geqslant 0.
  \end{equation}
We have already shown that $ 2 M K_{\max} < \dfrac{1}{2}$
in (\ref{Kmax}). This completes the proof of Lemma \ref{lemma-ind-ine}.
  
\vspace{0.2in}
\centerline{\bf \Large Appendix II. Proof of estimate (\ref{C1-integral})}

\vspace{0.2in}
In this appendix, we prove estimate (\ref{C1-integral}).
First, we state a useful inequality
\begin{equation}
|1 - e^{-2x}| \geqslant (1-e^{-2})|x|, \quad \mbox{for}\;\; 0 \leqslant x
\leqslant 1,
\label{ineq}
\end{equation}
which is a consequence of the fact that $(1-e^{-x})/x$ is a 
monotonely decreasing function for $x >0$.
We consider two cases.
If $\tau > 1$, we divide the integral into two parts as follows:
\begin{equation}
  \begin{array}{ll}
    & \int_0^{\tau} e^{ \gamma s} \left( 1 - e^{- 2 ( \tau - s )} 
    \right)^{- \frac{1+\gamma}{2}} d s\\[2mm]
    = & \int_0^{\tau - 1} + \int_{\tau - 1}^{\tau} e^{\gamma s}
    \left( 1 - e^{- 2( \tau - s )} \right)^{- \frac{1+\gamma}{2}} d s\\[2mm]
    \leqslant & \int_0^{\tau - 1} e^{\gamma s} \left( 1 - e^{- 2} 
    \right)^{- \frac{1+\gamma}{2}} d s + \int_{\tau - 1}^{\tau}
    e^{\gamma \tau}\left((1-e^{-2}) ( \tau - s ) \right)^{- \frac{1+\gamma}{2}} 
     d s\\[2mm]
    = & \left( 1 - e^{- 2} \right)^{- \frac{1+\gamma}{2}}
    \frac{e^{\gamma ( \tau - 1)} -1 }{\gamma} - e^{\gamma \tau}
    \left(1-e^{-2} \right)^{- \frac{1+\gamma}{2}} \left.
    \frac{2}{1-\gamma} \left( \tau - s \right)^{\frac{1-\gamma}{2}} \right|_{s=\tau -
    1}^{s=\tau}\\[2mm]
    = & \left( 1 - e^{- 2} \right)^{- \frac{1+\gamma}{2}}
    \frac{ e^{\gamma(\tau - 1)}-1}{\gamma} +\frac{2}{1-\gamma}
     e^{\gamma \tau }
     \left(1-e^{-2}\right)^{- \frac{1+\gamma}{2}}\\[2mm]
    \leqslant & \left ( \frac{2}{1-\gamma} + \frac{1}{\gamma} \right )
     \left(1-e^{-2}\right)^{- \frac{1+\gamma}{2}} e^{\gamma \tau},
  \end{array}
\label{C1}
\end{equation}
where we have used (\ref{ineq}). Thus we prove
\[
e^{-\tau} \int_0^{\tau} e^{ \gamma s} \left( 1 - e^{- 2 ( \tau - s )}
    \right)^{- \frac{1+\gamma}{2}} d s \leqslant c_1 e^{-(1-\gamma)\tau}
< c_1, \hspace{2em} \tmop{for} \;\; \tmop{all} \;\; \tau > 1,
\]
where $c_1$ is defined in (\ref{C1-def}).
For $\tau \leqslant 1$, we have by using (\ref{ineq})
\begin{eqnarray}
    e^{-\tau}\int_0^{\tau} e^{\gamma s} \left( 1 - e^{- 2 ( \tau - s )} 
     \right)^{- \frac{1+\gamma}{2}} d s & \leqslant & 
     e^{-(1-\gamma )\tau} 
   \int_0^{\tau} \left ( (1-e^{-2}) \left( \tau - s \right) \right)^{-
    \frac{1+\gamma}{2}} d s\nonumber\\[2mm]
   & \leqslant & \frac{2e^{-(1-\gamma)\tau}}{1-\gamma} \left(1-e^{-2}\right)^{- \frac{1+\gamma}{2}} 
    \tau^{\frac{1-\gamma}{2}} \nonumber\\[2mm] 
    & \leqslant & c_1 \tau^{\frac{1-\gamma}{2}} e^{-(1-\gamma)\tau} 
\leqslant   c_1 .    
\label{C1-appendix}
\end{eqnarray}
This proves (\ref{C1-integral}).

\vspace{0.2in}
\noindent
{\bf Acknowledgments.}
We would like to thank Prof. Congming Li for his comments and suggestions.
This work was in part supported by NSF under the NSF
FRG grant DMS-0353838 and ITR Grant ACI-0204932.

\bibliographystyle{amsplain}
\bibliography{bib}

\end{document}